# Empirical graph Laplacian approximation of Laplace–Beltrami operators: Large sample results

**Evarist Giné**[1,*] **and Vladimir Koltchinskii**[1,†]

*University of Connecticut and Georgia Institute of Technology*

**Abstract:** Let $M$ be a compact Riemannian submanifold of $\mathbf{R}^m$ of dimension $d$ and let $X_1, \ldots, X_n$ be a sample of i.i.d. points in $M$ with uniform distribution. We study the random operators

$$\Delta_{h_n, n} f(p) := \frac{1}{n h_n^{d+2}} \sum_{i=1}^{n} K(\frac{p - X_i}{h_n})(f(X_i) - f(p)), \ p \in M$$

where $K(u) := \frac{1}{(4\pi)^{d/2}} e^{-\|u\|^2/4}$ is the Gaussian kernel and $h_n \to 0$ as $n \to \infty$. Such operators can be viewed as graph laplacians (for a weighted graph with vertices at data points) and they have been used in the machine learning literature to approximate the Laplace-Beltrami operator of $M$, $\Delta_M f$ (divided by the Riemannian volume of the manifold). We prove several results on a.s. and distributional convergence of the deviations $\Delta_{h_n,n} f(p) - \frac{1}{|\mu|} \Delta_M f(p)$ for smooth functions $f$ both pointwise and uniformly in $f$ and $p$ (here $|\mu| = \mu(M)$ and $\mu$ is the Riemannian volume measure). In particular, we show that for any class $\mathcal{F}$ of three times differentiable functions on $M$ with uniformly bounded derivatives

$$\sup_{p \in M} \sup_{f \in \mathcal{F}} \left| \Delta_{h_n, p} f(p) - \frac{1}{|\mu|} \Delta_M f(p) \right| = O\left( \sqrt{\frac{\log(1/h_n)}{n h_n^{d+2}}} \right) \ \text{a.s.}$$

as soon as

$$n h_n^{d+2} / \log h_n^{-1} \to \infty \ \text{and} \ n h_n^{d+4} / \log h_n^{-1} \to 0,$$

and also prove asymptotic normality of $\Delta_{h_n,p} f(p) - \frac{1}{|\mu|} \Delta_M f(p)$ (functional CLT) for a fixed $p \in M$ and uniformly in $f$.

## 1. Introduction

Recently, there have been several developments in statistical analysis of data supported on a submanifold in a high dimensional space based on the idea of approximation of the Laplace-Beltrami operator of the manifold (and some more general operators that contain information not only about the geometry of the manifold, but also about the unknown density of data points) by empirical graph Laplacian operators. If $V$ is a finite set of vertices and $W := (w_{ij})_{i,j \in V}$ is a symmetric non-negative definite matrix of weights with $w_{ij} \geq 0$ ("adjacency matrix"), then the graph Laplacian of the weighted graph $(V, W)$ is defined as the matrix (operator) $L = D - W$, where $D$ is the diagonal matrix with the degrees of vertices

$$\deg(i) := \sum_{j \in V} w_{ij}, \ i \in V$$

on the diagonal. Such (unnormalized) graph Laplacians along with their normalized counterparts $\tilde{L} := I - D^{-1/2} L D^{-1/2}$ have been studied extensively in spectral graph theory. If now $X_1, \ldots, X_n$ are i.i.d. points uniformly distributed in a compact Riemannian submanifod $M$ of $\mathbf{R}^m$ of dimension $d < m$, it has been suggested in the literature to view $\{X_1, \ldots, X_n\}$ as the set $V$ of vertices of the graph and to define the weights as $w_{ij} \asymp e^{-\|X_i - X_j\|^2/4h^2}$ with a small parameter $h > 0$, to approximate the Laplace-Beltrami operator $\Delta_M$ of $M$, $\Delta_M(f) = \text{div}(\text{grad}(f))$. More precisely, the estimate is defined as

$$\Delta_{h_n, n} f(p) := \frac{1}{n h_n^{d+2}} \sum_{i=1}^n K\Big(\frac{p - X_i}{h_n}\Big)(f(X_i) - f(p)), \ p \in M$$

where $K(u) := \frac{1}{(4\pi)^{d/2}} e^{-\|u\|^2/4}$ is the Gaussian kernel and $h_n \to 0$ as $n \to \infty$ (if the functions $f$ are restricted to $V$, this can be viewed, up to a sign, as a graph Laplacian operator). We will call such operators *empirical graph Laplacians* and their limit as $n \to \infty$ on smooth functions $f$ is $\frac{1}{|\mu|} \Delta_M f(p)$, where $|\mu|$ is the Riemannian volume of $M$. There are numerous statistical applications of such an approximation of the manifold Laplacian by its empirical version. In particular, one can look at projections of the data on eigenspaces of the empirical Laplacian $\Delta_{h_n, n}$ (the technique sometimes called diffusion maps) in order to try to recover geometrically relevant features of the data (as in the method of spectral clustering) or use the kernels associated with this operator to approximate the heat kernel of the manifold and to use it to design kernel machines suitable, for instance, for classification of the manifold data.

[1]University of Connecticut, Department of Mathematics, U-3009, Storrs, CT 06269, USA, School of Mathematics, Georgia Inst. of Technology, Atlanta, GA 30332, USA, e-mail: gine@math.uconn.edu; vlad@math.gatech.edu
*Research partially supported by NSA Grant H98230-1-0075.
†Research partially supported by NSF Grant DMS-03-04861.

Convergence properties of empirical graph Laplacians have been first studied by Belkin and Niyogi [1] and Hein, Audibert and von Luxburg [8]. Our goal in this paper is to provide a more subtle probabilistic analysis of such operators. In particular, for proper classes of smooth functions $\mathcal{F}$ and for a fixed $p \in M$, we establish a functional CLT for $\sqrt{n h_n^{d+2}}(\Delta_{h_n, p} f(p) - \frac{1}{|\mu|} \Delta_M f(p))$, $f \in \mathcal{F}$, and also show that

$$\sup_{p \in M} \sup_{f \in \mathcal{F}} \Big|\Delta_{h_n, p} f(p) - \frac{1}{|\mu|} \Delta_M f(p)\Big| = O\Big(\sqrt{\frac{\log(1/h_n)}{n h_n^{d+2}}}\Big) \ \text{a.s.}$$

(under suitable assumptions on $h_n$). The asymptotic properties of empirical laplacians are closely related to the well developed theory of kernel density and kernel regression estimators, which can be viewed as examples of so called local empirical processes, as in [6]. Our proofs are essentially based on an extension of this type of results to the case of data on the manifolds (for kernel density estimation on manifolds, see, e.g., [11] and references therein). For simplicity, we are considering in the current paper only uniform distributions on manifolds and Gaussian kernels $K$, but more general types of operators that occur in the case when the distribution of the data is not uniform and more general kernels (as in the paper of Hein, Audibert and von Luxburg [8]) can be dealt with quite similarly using the methods of the paper.

## 2. Some geometric background

We refer to [4] for the basic definitions and notations from Riemannian geometry. Given a manifold $M$ and $p \in M$, $T_p(M)$ will denote the tangent space to $M$ at $p$, and $TM$ the tangent bundle. Let $M$ be a complete connected (embedded) Riemannian submanifold of $\mathbf{R}^m$, of dimension $d < m$, meaning that $M$ is a complete connected Riemannian manifold and that the inclusion map $\phi : M \mapsto \mathbf{R}^m$ is an isometric embedding, that is, (i) $\phi$ is differentiable and injective, (ii) $d\phi_p : T_p(M) \mapsto T_{\phi(p)}(\mathbf{R}^m)$ is an isometry onto its image, $T_{\phi(p)}(\phi(M))$, and (iii) $\phi$ is a homeomorphism onto $\phi(M)$ with the topology inherited from $\mathbf{R}^m$. When no confusion may arise, we identify $M$ with $\phi(M)$. $M$ being complete, by the Hopf and Rinow theorem (e.g., [4], p. 146) the closed bounded sets of $M$ are compact.

Given $p \in M$ and $v \in T_p(M)$, let $\gamma(t, p, v)$, $t > 0$, be the geodesic starting at $p$ with velocity $v$, $\gamma(0, p, v) = p$ and $\gamma'(0, p, v) = v$. The exponential map $\mathcal{E}_p : T_p(M) \mapsto M$ (the usual notation is $\exp_p$) is defined by $\mathcal{E}_p(v) = \gamma(1, p, v)$. This map is defined on all of $T_p(M)$ by the Hopf and Rinow theorem.

A normal neighborhood $V$ of $p \in M$ is one for which a) every point $q$ in $V$ can be joined to $p$ by a *unique* geodesic $\gamma(t, p, v)$, $0 \leq t \leq 1$, with $\gamma(0) = p$, $\gamma(1) = q$ and b) the exponential map centered at $p$, $\mathcal{E}_p$, is a diffeomorphism between a neighborhood of $0 \in T_p(M)$ and $V$. If $B \subset V$ is a normal ball of center $p$, that is, the image by the exponential map of a ball around zero in $T_p(M)$, then the unique geodesic joining $p$ to $q \in B$ is a minimizing geodesic, which means that if $d_M$ denotes the distance in $M$ and $|\cdot|$ denotes the norm of $T_p(M)$ defined by the Riemannian structure of $M$, then $d_M(p, \mathcal{E}_p(v)) = |v|$. Given an orthonormal basis $e_1, \ldots, e_d$ of $T_p(M)$, the normal coordinates centered at $p$ (or the $p$-normal coordinates) of $q \in V$ are the components $q_i^p = \langle \mathcal{E}_p^{-1}(q), e_i \rangle$ of $\mathcal{E}_p^{-1}(q)$ in this basis. (The super-index $p$ will be omitted when no confusion may arise, but we will need it when considering normal coordinates based at different points.) Every point in $M$ has a normal neighborhood. See [4], Propositions 2.7, and 3.6, pp. 64 and 70 for these facts. Actually, more is true (e.g., [4], Theorem 3.7 and Remark 3.8):

**Proposition 2.1.** *For every $p \in M$ there exist a neighborhood $W$ of $p$ and a number $\delta > 0$ such that:* (a) *for every $q$ in $W$, $\mathcal{E}_q$ is defined on the $\delta$ ball around $0 \in T_q(M)$, $B_\delta(0) \subset T_q(M)$, and $\mathcal{E}_q(B_\delta(0))$ is a normal ball for $q$,* (b) $W \subset \mathcal{E}_q(B_\delta(0))$ *($W$ is a normal neighborhood of all of its points), and* (c) *the function*

$$F(q, v) := (q, \mathcal{E}_q(v))$$

*is a diffeomorphism from $W_\delta := W \times B_\delta(0) = \{(q, v) \in TM : q \in W, |v| < \delta\}$ onto its image in $M \times M$ and $|dF|$ is bounded away from zero on $W_\delta$.*

Such a neighborhood $W$ of $p$ is called *totally or uniformly normal*. In particular, $\mathcal{E}_q(v)$ is jointly differentiable in $(q, v) \in W_\delta$ if $W$ is a uniformly normal neighborhood. Moreover, for every $q \in W$ and $v \in T_q(M)$ such that $|v| < \delta$, $d_M(q, \mathcal{E}_q(v)) = |v|$.

**Remark 2.1.** By shrinking $W$ and taking $\delta/2$ instead of $\delta$ if necessary, we can assume in Proposition 2.1 that the closure of $W$ and the closure of $W_\delta$ (which are compact because $M$ is complete) are contained in $W'$ and $W'_{\delta'}$ satisfying the properties described in the previous proposition. Moreover, we can also assume that for all $q$ in $W$, $\mathcal{E}_q(B_\delta(0))$ is contained in a strongly convex normal ball around $p$ (these points are at distances less than $2\delta$ from $p$, so this assumption can be met by further shrinking $W$ and taking a smaller $\delta$ if necessary, since every point in $M$ has a strongly convex geodesic ball, e.g. Proposition 4.2 in do Carmo, loc. cit.; strongly convex set: for any two points in the set, the minimizing geodesic joining them lies in the set). We will assume without loss of generality and without further mention that our uniformly normal neighborhoods $W$ satisfy these two conditions.

Let $W$ be a uniformly normal neighborhood of $p$ as in the remark, let $W'$ be a uniformly normal neighborhood of $p$ containing the closure of $W$, and let $0 < \delta < \delta'$ be as in the proposition and the remark. Let us choose an orthonormal basis $e_1, \ldots, e_d$ of $T_p(M)$ and define an orthonormal frame $e_1^q, \ldots, e_d^q$, $q \in W'$, by parallel transport of $e_1, \ldots, e_d$ from $p$ to $q$ along the unique minimizing geodesic joining $p$ and $q$. So, $e_1^q, \ldots, e_d^q$ is an orthonormal basis of $T_q(M)$ for each $q \in W'$. This frame depends differentiably on $q$ as parallel transport is differentiable (and preserves length and angle). So, we have on $W'$ a system of normal coordinates centered at $q$ for every $q \in W'$, namely, if $x \in \mathcal{E}_q(B_{\delta'}(0))$ is $x = \mathcal{E}_q(\sum_{i=1}^d v_i e_i^q)$, then the coordinates of $x$, $x_i^q$ are $x_i^q = v_i$, the components of $\mathcal{E}_q^{-1}(x)$. Let now $f$ be a differentiable function $f : M \mapsto \mathbf{R}$, and define $\tilde{f} : W'_{\delta'} \mapsto \mathbf{R}$ by

$$\tilde{f}(q,v) := f(\pi_2(F(q,v))) = f(\mathcal{E}_q(v)),$$

where $\pi_2$ is the projection of $M \times M$ onto its second component. This map is differentiable by the previous proposition. In particular, if we take as coordinates of $(q,v) \in W'_{\delta'}$ the normal coordinates centered at $p$ of $q$, $(q_1, \ldots, q_d) = (q_1^p, \ldots, q_d^p)$ and for $v \in T_q(M)$ the coordinates $v_1, \ldots v_d$ in the basis $e_1^q, \ldots, e_d^q$, which coincide with the normal coordinates centered at $q$ of $\mathcal{E}_q(v)$, then the real function of $2d$ variables (which we keep calling $\tilde{f}$; the same convention applies to other similar cases below)

$$\tilde{f}(q_1, \ldots, q_d, v_1, \ldots, v_d) = \tilde{f}(q,v)$$

is differentiable on the preimage of $W'_{\delta'}$ by this system of coordinates. Moreover, by compactness, each of its partial derivatives (of any order) is uniformly bounded on the preimage of $W_\delta$. If we denote by $x_i^q$ the normal coordinates centered at $q$, we obviously have that for each $r \in \mathbf{N}$ and $(i_1, \ldots, i_r) \in \{1, \ldots, d\}^r$,

$$\frac{\partial^r f}{\partial x_{i_1}^q \partial x_{i_2}^q \ldots \partial x_{i_r}^q}(x) = \frac{\partial^r \tilde{f}}{\partial v_{i_1} \partial v_{i_2} \ldots \partial v_{i_r}}(q_1, \ldots, q_d, x_1^q, \ldots, x_d^q).$$

We then conclude that

*each of the partial derivatives (any order) of $f$ with respect to the*

(2.1)

*$q$ − normal coordinates $x_i^q$ is uniformly bounded in $q \in W$ and $x \in \mathcal{E}_q(B_\delta(0))$.*

In particular, the error term in any limited Taylor development of $f$ in $q$-normal coordinates can be bounded uniformly in $q$ for all $|v| < \delta$, that is, if $P_k^q(x_1^q, \ldots, x_d^q)$ is the Taylor polynomial of degree $k$ in these coordinates, we have, for $q \in W$ and $|\mathcal{E}_q^{-1}(x)| < \delta$,

(2.2) $$|f(x) - P_k^q(x_1^q, \ldots, x_d^q)| \leq C_k(d_M(q,x))^{k+1},$$

where $C_k$ is a constant that depends only on $k$. Moreover, the coefficients of the polynomials $P_k^q$ are differentiable functions of $q$, in particular bounded on $W$. The $q$-uniformity of these Taylor developments for $|v| < \delta$ and the $q$-differentiability of their coefficients will be very useful below. We will apply these properties to the canonical (in $\mathbf{R}^m$) coordinates of the embedding $\phi$ and also to the functions $\langle \frac{\partial}{\partial x_i}, \frac{\partial}{\partial x_j}\rangle(x)$, where $x_i = x_i^p$ are the $p$-normal coordinates.

In what follows, we often deal with classes $\mathcal{F}$ of functions on $M$ whose partial derivatives up to a given order $k$ are uniformly bounded in $M$ or in a neighborhood $U$ of a point $p \in M$. In such cases, we say that $\mathcal{F}$ is *uniformly bounded up to the $k$-th order* in $M$ (or in $U$). Clearly, this property does not depend on the choice of normal (or even arbitrary) local coordinates. In the case when we choose an orthonormal frame $e_1^q, \ldots, e_d^q$ and define normal coordinates and the corresponding partial derivatives as described above, we can also deal with continuity of partial derivatives. We say that $\mathcal{F}$ is a *uniformly bounded and equicontinuous up to the $k$-th order* class of functions iff there exists a finite covering of $M$ with uniformly normal neighborhoods such that, in each neighborhood, the sets of partial derivatives of any order $\leq k$ are uniformly bounded and equicontinuous. This definition does not depend on the choice of orthonormal frames in the neighborhoods. Such classes are useful because the remainders in Taylor developments are uniform both in $q \in M$ and in $f \in \mathcal{F}$.

Consider now, for $q \in W'$ and $x \in \mathcal{E}_q(B_{\delta'}(0))$, the tangent vector fields $\frac{\partial}{\partial x_i^q}(x)$, $i = 1, \ldots, d$, and simply write $\frac{\partial}{\partial x_i}(x)$ for $\frac{\partial}{\partial x_i^p}(x)$. Taking the previous coordinates $q_i, v_j$ in $W'_{\delta'}$, denote by $\chi_i(\mathcal{E}_q(v)) = \chi_i(q_1, \ldots, q_d, v_1, \ldots, v_d)$ the $p$-normal coordinates of $\mathcal{E}_q(v)$, which are differentiable. By the chain rule,

$$\frac{\partial}{\partial x_i^q}(x) = \sum_{j=1}^{d} \frac{\partial \chi_j}{\partial v_i}(q, \mathcal{E}_q^{-1}x) \frac{\partial}{\partial x_j}(x).$$

Hence, if $g_{ij}^q(x)$ are the components of the metric tensor at $x$ in $q$-normal coordinates, we have

$$g_{ij}^q(x) = \langle \frac{\partial}{\partial x_i^q}, \frac{\partial}{\partial x_j^q}\rangle(x) = \sum_{1 \leq r,s \leq d} \frac{\partial \chi_r}{\partial v_i}(q, \mathcal{E}_q^{-1}x)\frac{\partial \chi_s}{\partial v_j}(q, \mathcal{E}_q^{-1}x)\langle \frac{\partial}{\partial x_i}, \frac{\partial}{\partial x_j}\rangle(x).$$

By (2.2), we conclude that if $P_k^q(x_1^q, \ldots, x_d^q)$ is the Taylor polynomial of degree $k$ in the expansion of $g_{ij}^q(x)$ in $q$-normal coordinates, then there are constants $C_k$ that depend only on $k$ such that, for all $q \in W$ and $x \in \mathcal{E}_q(B_\delta(0))$,

(2.3) $$|g_{ij}^q(x) - P_k^q(x_1^q, \ldots, x_d^q)| \leq C_k(d_M(q,x))^{k+1}.$$

This will also be useful below. These remarks allow us to strengthen several results based on Taylor expansions by making them uniform in $q$, as follows.

**Proposition 2.2.** *Given $p \in M$, let $W$ and $W_\delta$ be as in Remark 2.1, and consider for each $q \in W$, the $q$-normal system of coordinates defined above. Then,*
*(a) for every $q \in W$ the components $g_{ij}^q(x_1^q, \ldots, x_d^q)$ of the metric tensor in $q$-normal coordinates admit the following expansion, uniform in $q$ and in $x \in \mathcal{E}_q(B_\delta(0))$ ($B_\delta(0) \in T_q(M)$):*

(2.4) $$g_{ij}^q(x_1^q, \ldots, x_d^q) = \delta_{ij} - \frac{1}{3}R_{irsj}^q(0)x_r^q x_s^q + O(d_M^3(q,x)),$$

*(Einstein notation) where $R_{irsj}^q(0)$ are the components of the curvature tensor at $q$ in $q$-normal coordinates, and, as a consequence, the following expansion of the volume element is also uniform in $q$ and $x$:*

(2.5) $$\sqrt{\det(g_{ij}^q)}(x_1^q, \ldots, x_d^q) = 1 - \frac{1}{6}\operatorname{Ric}_{rs}^q(0)x_r^q x_s^q + O(d_M^3(q,x)),$$

where $\mathrm{Ric}^q_{rs}(0)$ are the components of the Ricci tensor at $q$ in $q$-normal coordinates.
(b) *There exists $C < \infty$ such that for all $q \in W$ and $x \in \mathcal{E}_q(B_\delta(0))$,*

$$(2.6) \qquad 0 \le d_M^2(q,x) - \|\phi(q) - \phi(x)\|^2 \le C d_M^4(q,x).$$

(c) *For each $1 \le \alpha \le m$, the $\alpha$-th component in canonical coordinates of $\mathbf{R}^m$ of $\phi(\mathcal{E}_q(v))$, $\phi_\alpha(\mathcal{E}_q(v))$, admits the following expansion in $q$-normal coordinates $v_i$ of $\mathcal{E}_q(v)$, uniform in $q \in W$ and $|v| < \delta$,*

$$(2.7) \qquad \phi_\alpha(\mathcal{E}_q(v)) - \phi_\alpha(q) = \frac{\partial \tilde{\phi}_\alpha}{\partial v_i}(q,0) v_i + \frac{1}{2} \frac{\partial^2 \tilde{\phi}_\alpha}{\partial v_i \partial v_j}(q,0) v_i v_j + O(|v|^3),$$

$\alpha = 1, \ldots, m$, where $\tilde{\phi}(q,v) = \phi(\mathcal{E}_q(v))$.

Note that $\sum_{i=1}^d \frac{\partial \tilde{\phi}_\alpha}{\partial v_i}(q,0) v_i$ are the $\mathbf{R}^m$-canonical coordinates centered at $\phi(q)$ of the vector $d\phi_q(v) \in T_{\phi(q)}(\phi(M)) \subset \mathbf{R}^m$ since $g^q_{ij}(q) = \delta_{ij}$. Hence, if we identify the tangent space to $\phi(M)$ at $\phi(q)$ with an affine subspace of $\mathbf{R}^m$, part c) says that the difference between $\phi(\mathcal{E}_q(v)) \in \mathbf{R}^m$ and the tangent vector to the geodesic $\phi(\gamma(t,q,v))$ at $\phi(q)$ ($t = 0$), $\phi(q) + d\phi_q(v)$, is a vector of the form

$$\Big( \frac{1}{2} \sum_{1 \le i,j \le d} \frac{\partial^2 \tilde{\phi}_\alpha}{\partial v_i \partial v_j}(q,0) v_i v_j : \alpha = 1, \ldots, m \Big) + O(|v|^3)$$

where $O(|v|^3)$ is uniform in $|v| < \delta$ and $q$. Ignoring the embedding, this gives an expansion of the exponential map as

$$(2.7') \qquad \mathcal{E}_q(v) = q + v + Q_q(v,v) + O(|v|^3)$$

uniform in $q \in W$ and $|v| < \delta$, where $Q_q$ is a $\mathbf{R}^m$-valued bilinear map on $T_q(M)$ (actually, on $T_{\phi(q)}(\phi(M))$) that depends differentiably on $q$, hence uniformly bounded in $q \in W$.

*Proof of Proposition 2.2.* (a) follows from the expansions of $g^q_{ij}$ and $\sqrt{\det(g^q_{ij})}$ in $q$-normal coordinates (e.g. in [12], p. 41), the expansion of its determinant (e.g., [12], p. 45), and the uniformity provided by (2.3).

(c) follows by direct application of (2.2) to $f = \phi_\alpha$, $\alpha = 1, \ldots, m$.

(b) Following Smolyanov, Weizsäcker and Wittich [13], for $q \in W$ and $x = \mathcal{E}_q(v)$, $|v| < \delta$, and applying (2.2) for $f = \phi_\alpha$, $\alpha = 1, \ldots, m$, we have

$$0 \le \frac{d_M^2(q,x) - \|\phi(q) - \phi(x)\|^2}{d_M^4(q,x)} = \frac{|v|^2 - \sum_{\alpha=1}^m \big(\phi_\alpha(\mathcal{E}_q(v)) - \phi_\alpha(q)\big)^2}{|v|^4}$$

$$= \frac{|v|^2 - \sum_{\alpha=1}^m \big(\frac{\partial \tilde{\phi}_\alpha}{\partial v_i}(q,0) v_i + \frac{1}{2} \frac{\partial^2 \tilde{\phi}_\alpha}{\partial v_i \partial v_j}(q,0) v_i v_j + \frac{1}{6} \frac{\partial^3 \tilde{\phi}_\alpha}{\partial v_i \partial v_j \partial v_k}(q,0) v_i v_j v_k\big)^2}{|v|^4}$$

$$+ O(|v|),$$

where the term $O(|v|)$ is dominated by $C_4 |v|$ for a constant $C_4$ that does not depend on $q$ or $v$. But now, continuing the proof in this reference, which consists in developing and simplifying the ratio above, we obtain that, uniformly in $q \in W$, $x = \mathcal{E}_q(v)$, $|v| < \delta$,

$$0 \le \frac{d_M^2(q,x) - \|\phi(q) - \phi(x)\|^2}{d_M^4(q,x)} = \frac{1}{12} \sum_\alpha \frac{\big(\frac{\partial^2 \tilde{\phi}_\alpha}{\partial v_i \partial v_j}(q,0) v_i v_j\big)^2}{|v|^4} + O(|v|),$$

and note also that, by compactness, the main term is bounded by a fixed finite constant in this domain. □

Although we have been using [4] as our main reference on Riemannian geometry, another nice user-friendly reference for the exponential map in particular and Riemannian manifolds in general is [9]. We thank Jesse Ratzkin for reading this section and making comments (of course, any mistakes are ours).

## 3. Approximation of the Laplacian by averaging kernel operators

Let $M$ be a compact connected Riemannian submanifold of $\mathbf{R}^m$, $m > d$ (if $M$ is compact, it is automatically embedded, that is, conditions (i) and (ii) on the immersion $\phi$ imply that $\phi$ is a homeomorphism onto its image). [ See a remark at the end of this section for a relaxation of this condition.] Let $\mu$ be its Riemannian volume measure and $|\mu| = \mu(M)$. Let $K : \mathbf{R}^m \mapsto \mathbf{R}$ be the Gaussian kernel of $\mathbf{R}^m$,

$$(3.1) \qquad K(x) = \frac{1}{(4\pi)^{d/2}} e^{-\|x\|^2/4},$$

where $\|x\|$ is the norm of $x$ in $\mathbf{R}^m$. Let $X$ be a random variable taking values in $M$ with law the normalized volume element, $\mu/|\mu|$, and let $f : M \mapsto \mathbf{R}$ be a differentiable function. The object of this section is to show that the Laplace-Beltrami operator or Laplacian of $M$,

$$\Delta_M f(p) = \operatorname{div} \operatorname{grad}(f)(p)$$

(in coordinates, $\Delta_M(f) = \frac{1}{\sqrt{\det(g_{ij})}} \frac{\partial}{\partial x_i}(g^{ij}\sqrt{\det(g_{ij})}\frac{\partial f}{\partial x_j})$, where $(g^{ij}) = (g_{ij})^{-1}$) can be approximated, uniformly in $f$ (with some partial derivatives bounded), and in $p \in M$, by the *averaging kernel operator*

$$(3.2) \qquad \Delta_{h_n} f(p) := \frac{1}{h_n^{d+2}} E\left[K\left(\frac{\phi(p) - \phi(X)}{h_n}\right)(f(X) - f(p))\right]$$

with rates depending on $h_n \to 0$. Note that, by the expansion (2.4) of the metric tensor in normal coordinates centered at $p$, we have, in these coordinates,

$$(3.3) \qquad \Delta_M f(p) = \sum_{i=1}^d \frac{\partial^2 f}{\partial x_i^2}(p).$$

(where $p = (0, \ldots, 0)$ in these coordinates).

With some abuse of notation, given $p \in M$, we denote the derivatives with respect to the components of $v$ of $\tilde{f}(p,v) = f \circ \mathcal{E}_p(v)$ at $(p,v)$, $v = \mathcal{E}_p^{-1}(x)$, by $f'(x)$, $f''(x)$, etc. (so, for instance, if $x = \mathcal{E}_p(v)$, $f'(x) = (\frac{\partial \tilde{f}}{\partial v_1}(p,v), \ldots, \frac{\partial \tilde{f}}{\partial v_d}(p,v)))$ (in fact, $f^{(k)}(x)$ depends on $p$ and therefore it should have been denoted $f_p^{(k)}(x)$, but in the context we are using this notation $p$ is typically fixed, so, we will drop $p$, hopefully, without causing a confusion).

**Theorem 3.1.** *We have, for any $p$, any normal neighborhood $U_p$ of $p$ and a class $\mathcal{F}$ uniformly bounded up to the third order in $U_p$, that*

$$(3.4) \qquad \sup_{f \in \mathcal{F}}\left|\Delta_{h_n} f(p) - \frac{1}{|\mu|}\Delta_M f(p)\right| = O(h_n).$$

as $h_n \to 0$. Moreover, for any class of functions uniformly bounded up to the third order in $M$,

$$(3.5) \qquad \sup_{f \in \mathcal{F}} \sup_{p \in M} \left| \Delta_{h_n} f(p) - \frac{1}{|\mu|} \Delta_M f(p) \right| = O(h_n).$$

as $h_n \to 0$.

*Proof.* $M$ being regular, the embedding $\phi$ is a homeomorphism of $M$ onto $\phi(M)$, and $M$ being compact, the uniformities defined respectively on $M$ by the intrinsic metric $d_M(p,q)$ and by the metric from $\mathbf{R}^m$, $d_{\mathbf{R}^m}(p,q) := \|\phi(p) - \phi(q)\|$ coincide (e.g., Bourbaki (1940), Theorem II.4.1, p. 107), that is, given $\varepsilon > 0$ there exists $\delta > 0$ such that if $d_M(p,q) < \delta$ for $p, q \in M$, then $d_{\mathbf{R}^m}(p,q) < \varepsilon$, and conversely. Hence, in Proposition 2.2, we can replace $B_\delta(0) \subset T_q(M)$ by $B'_{\delta'}(0) := \mathcal{E}_q^{-1}\{x \in M : \|\phi(q) - \phi(x)\| < \delta'\}$ for some $\delta'$ depending on $\delta$ but not on $p$ or $q$. From here on, we identify $M$ with $\phi(M)$ (that is, we leave $\phi$ implicit). Let $p \in M$. Given $h_n \searrow 0$, let

$$(3.6) \qquad \mathcal{B}_n := \{x \in M : \|p - x\| < Lh_n(\log h_n^{-1})^{1/2}\}$$

for a constant $L$ to be chosen later. As soon as $Lh_n(\log h_n^{-1})^{1/2} < \delta'$, the neighborhood of $0 \in T_p(M)$,

$$\tilde{\mathcal{B}}_n := \mathcal{E}_p^{-1} \mathcal{B}_n$$

is well defined, and, by (2.6), since $|v| = d_M(p, \mathcal{E}_p(v))$, we have on $\tilde{\mathcal{B}}_n$ that

$$|v|^2 \geq \|p - \mathcal{E}_p(v)\|^2 \geq |v|^2(1 - C|v|^2)$$

with $C$ independent of $p \in M$. Hence, for all $n \geq N_0$, for some $N_0 < \infty$ independent of $p$, we have

$$(3.7) \qquad \begin{aligned} \{v \in T_p M : |v| < Lh_n(\log h_n^{-1})^{1/2}\} &\subseteq \tilde{\mathcal{B}}_n \\ &\subseteq \{v \in T_p M : |v| < 2Lh_n(\log h_n^{-1})^{1/2}\}, \end{aligned}$$

where the coefficient 2 can be replaced by $\lambda_n \to 1$. Assume $n \geq N_0$.

By the definitions of $K$ and $\mathcal{B}_n$,

$$(3.8) \qquad \begin{aligned} E\left| K\left(\frac{p-X}{h_n}\right)(f(X) - f(p))I(X \in M \setminus \mathcal{B}_n) \right| \\ \leq \frac{2\|f\|_\infty}{(4\pi)^{d/2}} \int_{M \setminus \mathcal{B}_n} e^{-\|p-x\|^2/4h_n^2} \frac{d\mu(x)}{|\mu|} \\ \leq \frac{2\|f\|_\infty}{(4\pi)^{d/2}} h_n^{L^2/4}. \end{aligned}$$

Taking into account that the measure $\mu$ has density $\sqrt{\det(g_{ij})}$ in $p$-normal coordinates (hence on $\mathcal{B}_n$), we have

$$(3.9) \qquad \begin{aligned} E\left[ K\left(\frac{p-X}{h_n}\right)(f(X) - f(p))I(X \in \mathcal{B}_n) \right] \\ = \frac{1}{(4\pi)^{d/2}|\mu|} \int_{\tilde{\mathcal{B}}_n} e^{-\|p-\mathcal{E}_p(v)\|^2/4h_n^2} (f(\mathcal{E}_p(v)) - f(\mathcal{E}_p(0))) \sqrt{\det(g_{ij})}(v) dv. \end{aligned}$$

With the notation introduced just before the statement of the theorem, the Taylor expansion of $f$ in $p$-normal coordinates can be written as

$$f(\mathcal{E}_p(v)) - f(\mathcal{E}_p(0)) = \langle f'(p), v \rangle + \frac{1}{2} f''(p)(v,v) + \frac{1}{3!} f'''(\xi_v)(v,v,v).$$

where $\xi_v = \mathcal{E}_p(\theta_v v)$ for some $\theta_v \in [0, 1]$. Next we will estimate the three terms that result from combining this Taylor development with equation (3.9). Recall that, by Proposition 2.2, there are $C_1$ and $C$ independent of $p$ such that

$$(3.10) \quad \sqrt{\det(g_{ij})}(v) \leq 1 + C_1|v|^2, \quad \frac{1}{2}|v|^2 \leq |v|^2 - C|v|^4 \leq \|p - \mathcal{E}_p(v)\|^2 \leq |v|^2$$

for $v \in \tilde{\mathcal{B}}_n$, and recall also (3.7) on the size of $\tilde{\mathcal{B}}_n$. Using these facts and the development of the exponential about $-|v|^2/4h_n^2$ immediately gives

$$\int_{\tilde{\mathcal{B}}_n} e^{-\|p-\mathcal{E}_p(v)\|^2/4h_n^2} \langle f'(p), v\rangle \sqrt{\det(g_{ij})}(v) dv = \int_{\tilde{\mathcal{B}}_n} e^{-|v|^2/4h_n^2} \langle f'(p), v\rangle dv + R_n$$

where

$$|R_n| \leq \int_{\tilde{\mathcal{B}}_n} \left(e^{-(|v|^2 - C|v|^4)/4h_n^2} - e^{-|v|^2/4h_n^2}\right)|f'(p)||v|dv$$
$$+ C_1 \int_{\tilde{\mathcal{B}}_n} e^{-|v|^2/8h_n^2}|f'(p)||v|^3 dv$$
$$\leq \int_{\tilde{\mathcal{B}}_n} e^{-|v|^2/8h_n^2}|f'(p)|(C|v|^5/(4h_n^2) + C_1|v|^3) dv$$
$$\leq h_n^{3+d} \int_{\mathbf{R}^d} e^{-|v|^2/8}|f'(p)|(C|v|^5/4 + C_1|v|^3) dv$$
$$= D|f'(p)|h_n^{3+d},$$

and $D$ only depends on $C$, $C_1$ and $d$. Moreover, since $\mathcal{B}_n^c \subseteq \{|v| \geq Lh_n(\log h_n^{-1})^{1/2}\}$ and

$$\int_{\mathbf{R}^d} e^{-|v|^2/4h_n^2} \langle f'(p), v\rangle dv = 0,$$

we also have

$$\left|\int_{\tilde{\mathcal{B}}_n} e^{-|v|^2/4h_n^2} \langle f'(p), v\rangle dv\right| = \left|\int_{\tilde{\mathcal{B}}_n^c} e^{-|v|^2/4h_n^2} \langle f'(p), v\rangle dv\right|$$
$$\leq |f'(p)| \int_{|v| \geq Lh_n(\log h_n^{-1})^{1/2}} e^{-|v|^2/4h_n^2}|v| dv$$
$$= |f'(p)| h_n^{1+d} \int_{|u| \geq L(\log h_n^{-1})^{1/2}} e^{-|u|^2/4}|u| du$$
$$= C_d|f'(p)| h_n^{1+d} \int_{r \geq L(\log h_n^{-1})^{1/2}} e^{-r^2/4} r^d dr$$
$$\leq C_d'|f'(p)| L^{d-1} h_n^{1+d+L^2/4} (\log h_n^{-1})^{(d-1)/2},$$

where $C_d$ and $C_d'$ are constants depending only on $d$. Collecting terms and assuming

$$L^2/4 > 2,$$

we obtain

$$(3.11) \quad \left|\int_{\tilde{\mathcal{B}}_n} e^{-\|p-\mathcal{E}_p(v)\|^2/4h_n^2} \langle f'(p), v\rangle \sqrt{\det(g_{ij})}(v) dv\right| \leq D_2|f'(p)|h_n^{3+d},$$

for all $n \geq N_0$, and where $D_2$ does not depend on $p$.

The remainder term is of a similar order if $|f'''|$ is uniformly bounded in a neighborhood of $p$: if $c$ is such a bound,

$$\left| \int_{\tilde{\mathcal{B}}_n} e^{-\|p-\mathcal{E}_p(v)\|^2/4h_n^2} f'''(\xi_v)(v,v,v)\sqrt{\det(g_{ij})}(v)dv \right| \tag{3.12}$$
$$\leq c \int_{\tilde{\mathcal{B}}_n} e^{-|v|^2/8h_n^2} |v|^3 |1 + C_1 h_n^2 \log h_n^{-1}| dv \leq D_3 c h_n^{3+d},$$

where $D_3$ does not depend on $f$ or $p$ (as long as $n \geq N_0$).

Finally, we consider the second term, which is the one that gives the key relationship to the Laplacian. Proceeding as we did for the first term, we see that

$$\int_{\tilde{\mathcal{B}}_n} e^{-\|p-\mathcal{E}_p(v)\|^2/4h_n^2} f''(p)(v,v)\sqrt{\det(g_{ij})}(v)dv \tag{3.13}$$
$$= h_n^{d+2} \int_{\mathbf{R}^d} e^{-|v|^2/4} f''(p)(v,v)dv + R_n,$$

where now

$$(3.14) \quad |R_n| \leq D_4 |f''(p)| h_n^{4+d} + D_5 |f''(p)| h_n^{2+d+L^2/4} (\log h_n^{-1})^{d/2} \leq D_6 |f''(p)| h_n^{4+d}$$

if

$$L^2/4 > 2$$

and $n \geq N_0$, where the constants $D$ do not depend on $f$ or $p$. Now, by definition

$$f''(p) = (f \circ \mathcal{E}_p)''(0) = \left( \frac{\partial^2(f \circ \mathcal{E}_p)}{\partial v_i \partial v_j}(0) \right)_{i,j=1}^d = \left( \frac{\partial^2 f}{\partial x_i \partial x_j}(p) \right)_{i,j=1}^d,$$

so that, on account of (3.3),

$$\int_{\mathbf{R}^d} e^{-|v|^2/4} f''(p)(v,v)dv = \left( \int_{\mathbf{R}^d} e^{-|v|^2/4} v_1^2 dv \right) \sum_{i=1}^d \frac{\partial^2 f}{\partial x_i^2}(p) \tag{3.15}$$
$$= 2(4\pi)^{d/2} \sum_{i=1}^d \frac{\partial^2 f}{\partial x_i^2}(p) = 2(4\pi)^{d/2} \Delta_M f(p).$$

Combining the bounds (3.11), (3.12), (3.13)-(3.14) and the identity (3.15) with (3.9), we obtain the first part of the theorem. Note that we need to choose $L$ such that $L^2/4 > 2$ and then $N_0$ such that $Lh_{N_0}(\log h_{N_0}^{-1})^{1/2} < \delta$, and that with these choices the bounds obtained on the terms that tend to zero in the proof depend only on the sup of the derivatives of $f$ and on the sup of certain differentiable functions of $(q,v)$ on $W_\delta$, where $W$ is a uniformly normal neighborhood of $p$ and $\delta$ the corresponding number from Proposition 2.2 and Remark 2.1. These bounds are the same if we replace $p$ by any $q \in W$ by Proposition 2.2. $M$ being compact, it can be covered by a finite number of uniformly normal neighborhoods $W_i$, $i \leq k$, with numbers $\delta_i$ as prescribed in Proposition 2.1 and Remark 2.1. Taking $\delta$ to be the minimum of $\delta_i$, $i = 1, \ldots, k$, and the constants in the bounds in the first part of the proof as the maximum of the constants in these bounds for each of the $k$ neighborhoods, the above estimates work uniformly on $q \in M$, giving the second part of the theorem. □

**Remark 3.1.** (1) Obviously, the first part of the theorem, namely the limit (3.4), does not require the manifold $M$ to be compact. (2) If instead of assuming existence and boundedness of the third order partial derivatives in a neighborhood of $p$ we assume that the second order derivatives are continuous in a neighborhood of $p$, then we can proceed as in the above proof except for the remainder term (3.12), that now can be replaced by

$$(3.16) \quad \left| \int_{\tilde{\mathcal{B}}_n} e^{-\|p - \mathcal{E}_p(v)\|^2 / 4h_n^2} (f''(\xi_v) - f''(p))(v,v) \sqrt{\det(g_{ij})}(v) dv \right| = o(h_n^{2+d}).$$

Hence, in this case we still have

$$(3.17) \quad \Delta_{h_n} f(p) \to \frac{1}{|\mu|} \Delta_M f(p) \text{ as } h_n \to 0.$$

A similar observation can be made regarding (3.5).

**Remark 3.2.** Suppose $N$ is a compact Riemannian $d$-dimensional submanifold of $\mathbf{R}^m$ *with boundary* (for the definition, see [12], p. 70-71). The Riemannian volume measure $\mu$ is still finite. Then, Theorem 3.1 is still true if $X$ is a $N$-valued random variable with law $\mu/|\mu|$ with $|\mu| = \mu(N)$, and $M$ a compact subset of $N$ interior to $N$. The proof is essentially the same.

The first part of Theorem 3.1, without uniformity in $f$, is proved in a more general setting in [8].

Theorem 3.1 provides the basis for the estimation of the Laplacian of $M$ by independent sampling from the space according to the normalized volume element, which is what we do for the rest of this article.

## 4. Pointwise approximation of the Laplacian by graph Laplacians

Let $M$ be a compact Riemannian submanifold of $\mathbf{R}^d$ (or, in more generality, let $M$ be as in Remark 3.2), and let $X, X_i, i \in \mathbf{N}$, be independent identically distributed random variables with law $\mu/|\mu|$. The 'empirical counterpart' of the averaging kernel operator from Section 3 corresponding to such a sequence is the so called *graph Laplacian*

$$(4.1) \quad \Delta_{h_n, n} f(p) := \frac{1}{nh_n^{d+2}} \sum_{i=1}^n K\left(\frac{p - X_i}{h_n}\right)(f(X_i) - f(p)),$$

with $K$ given by (3.1) (other kernels are possible).

We begin with the pointwise central limit theorem for a single function $f$, as a lemma for the CLT uniform in $f$.

**Proposition 4.1.** *Assume $f$ has partial derivatives up to the third order continuous in a neighborhood of $p$. Let $h_n \to 0$ be such that $nh_n^d \to \infty$ and $nh_n^{d+4} \to 0$. Then,*

$$(4.2) \quad \sqrt{nh_n^{d+2}} \left[\Delta_{h_n, n} f(p) - \frac{1}{|\mu|} \Delta_M f(p)\right] \to sg \text{ in distribution,}$$

*where $g$ is a standard normal random variable and*

$$(4.3) \quad s^2 = \frac{1}{(4\pi)^d |\mu|} \int_{\mathbf{R}^d} e^{-|v|^2/2} \left(\sum_{j=1}^d \frac{\partial f}{\partial x_j^p}(p) v_j\right)^2 dv = \frac{1}{2^d (2\pi)^{d/2} |\mu|} \sum_{j=1}^d \left|\frac{\partial f}{\partial x_j^p}(p)\right|^2.$$

*Proof.* Since by Theorem 3.1,

$$\sqrt{nh_n^{d+2}}\Big(\Delta_M f(p) - \Delta_{h_n} f(p)\Big) = O\Big(\sqrt{nh_n^{d+4}}\Big) \to 0,$$

it suffices to prove that the sequence

(4.4)
$$Z_n = \sqrt{nh_n^{d+2}} \Bigg[ \frac{1}{nh_n^{d+2}} \sum_{i=1}^n \Bigg( K\Big(\frac{p-X_i}{h_n}\Big)(f(X_i) - f(p)) \\ - EK\Big(\frac{p-X}{h_n}\Big)(f(X_i) - f(p))\Bigg) \Bigg]$$

is asymptotically centered normal with variance $s^2$. To prove this we first observe that we can restrict to $X_i \in \mathcal{B}_n$ because, as in (3.8),

(4.5)
$$\frac{1}{h_n^{d+2}} E\Bigg[K^2\Big(\frac{p-X_i}{h_n}\Big)(f(X_i) - f(p))^2 I(X \in M \setminus \mathcal{B}_n)\Bigg] \\ \le \frac{4\|f\|_\infty^2}{(4\pi)^d h_n^{d+2}} \int_M e^{-L^2(\log h_n^{-1})/2} d\mu/|\mu| = \frac{4\|f\|_\infty^2}{(4\pi)^d} h_n^{L^2/2 - (d+2)} \to 0$$

if we take $L^2/2 > d+2$. Now, on the restriction to $\mathcal{B}_n$ we replace $f(X_i) - f(p)$ by its Taylor expansion up to the second order plus remainder, as done in the proof of Theorem 3.1. The second term and the remainder parts, namely

(4.6)
$$Z_{n,2} := \sqrt{nh_n^{d+2}} \Bigg[ \frac{1}{nh_n^{d+2}} \sum_{i=1}^n \Bigg( K\Big(\frac{p-X_i}{h_n}\Big) f''(p)(\mathcal{E}_p^{-1}(X_i), \mathcal{E}_p^{-1}(X_i)) I(X_i \in \mathcal{B}_n) \\ - EK\Big(\frac{p-X}{h_n}\Big) f''(p)(\mathcal{E}_p^{-1}(X), \mathcal{E}_p^{-1}(X)) I(X \in \mathcal{B}_n)\Bigg) \Bigg]$$

and

(4.7)
$$Z_{n,3} := \sqrt{nh_n^{d+2}} \Bigg[ \frac{1}{nh_n^{d+2}} \sum_{i=1}^n \Bigg( K\Big(\frac{p-X_i}{h_n}\Big) f'''(\xi_i) \\ \times (\mathcal{E}_p^{-1}(X_i), \mathcal{E}_p^{-1}(X_i), \mathcal{E}_p^{-1}(X_i)) I(X_i \in \mathcal{B}_n) \\ - EK\Big(\frac{p-X}{h_n}\Big) f'''(\xi)(\mathcal{E}_p^{-1}(X), \mathcal{E}_p^{-1}(X), \mathcal{E}_p^{-1}(X)) I(X \in \mathcal{B}_n)\Bigg) \Bigg]$$

tend to zero in probability: the estimates (3.10) give

(4.8)
$$EZ_{n,2}^2 \le \frac{1}{(4\pi)^d |\mu| h_n^{d+2}} \int_{\tilde{\mathcal{B}}_n} e^{-\|p-\mathcal{E}_p(v)\|^2/2h_n^2} (f''(p)(v,v))^2 \sqrt{\det(g_{ij})}(v) dv \\ \le \frac{2}{(4\pi)^d |\mu| h_n^{d+2}} \int_{\tilde{\mathcal{B}}_n} e^{-|v|^2/3h_n^2} |f''(p)|^2 |v|^4 dv \\ \le \frac{2h_n^{d+4}}{(4\pi)^d |\mu| h_n^{d+2}} \int_{\mathbf{R}^d} e^{-|v|^2/3} |f''(p)|^2 |v|^4 dv = O(h_n^2) \to 0,$$

and, with $c = \sup_{x \in U} |f'''(x)|$,

$$
\begin{aligned}
E|Z_{n,3}| &\leq \frac{2cn}{(4\pi)^{d/2}|\mu|\sqrt{nh_n^{d+2}}} \int_{\tilde{\mathcal{B}}_n} e^{-\|p - \mathcal{E}_p(v)\|^2/4h_n^2} |v|^3 \sqrt{\det(g_{ij})}(v) dv \\
&\leq \frac{3cn}{(4\pi)^{d/2}|\mu|\sqrt{nh_n^{d+2}}} \int_{\tilde{\mathcal{B}}_n} e^{-|v|^2/5h_n^2} |v|^3 dv \\
&\leq \frac{3h_n^{3+d}cn}{(4\pi)^{d/2}|\mu|\sqrt{nh_n^{d+2}}} \int_{R^d} e^{-|v|^2/5} |v|^3 dv \\
&= O\left(\sqrt{nh_n^{d+4}}\right) \to 0.
\end{aligned}
\tag{4.9}
$$

Finally, we show that the linear term part,

$$
\begin{aligned}
Z_{n,1} := \sqrt{nh_n^{d+2}} \Bigg[ \frac{1}{nh_n^{d+2}} \sum_{i=1}^{n} \Bigg( K\left(\frac{p - X_i}{h_n}\right) \langle f'(p), \mathcal{E}_p^{-1}(X_i) \rangle I(X_i \in \mathcal{B}_n) \\
- EK\left(\frac{p - X}{h_n}\right) \langle f'(p), \mathcal{E}_p^{-1}(X) \rangle I(X \in \mathcal{B}_n) \Bigg) \Bigg]
\end{aligned}
\tag{4.10}
$$

is asymptotically $N(0, s^2)$. Since, by (3.11),

$$
\begin{aligned}
\sqrt{\frac{n}{h_n^{d+2}}} EK\left(\frac{p - X}{h_n}\right) \langle f'(p), \mathcal{E}_p^{-1}(X) \rangle I(X \in \mathcal{B}_n) \\
= \sqrt{\frac{n}{h_n^{d+2}}} O(h_n^{d+3}) = O\left(\sqrt{nh_n^{4+d}}\right) \to 0,
\end{aligned}
$$

and since, by computations similar to the ones leading to (3.11),

$$
\begin{aligned}
\frac{1}{h_n^{d+2}} E \left[ K\left(\frac{p - X}{h_n}\right) \langle f'(p), \mathcal{E}_p^{-1}(X) \rangle I(X \in \mathcal{B}_n) \right]^2 \\
= \frac{1}{(4\pi)^d |\mu| h_n^{d+2}} \int_{\tilde{\mathcal{B}}_n} e^{-\|p - \mathcal{E}_p(v)\|^2/2h_n^2} \langle f'(p), v \rangle^2 \sqrt{\det(g_{ij})}(v) dv \\
= \frac{1}{(4\pi)^d |\mu|} \int_{\mathbf{R}^d} e^{-|v|^2/2} \langle f'(p), v \rangle^2 dv \\
+ \frac{1}{h_n^{d+2}} \left( O\left(h_n^{d+4}\right) + O\left(h_n^{2+d+L^2/4}(\log h_n^{-1})^{(d+1)/2}\right) \right),
\end{aligned}
$$

we have that, taking $L^2/4 > 2 + d$,

$$
\lim_{n \to \infty} EZ_{n,1}^2 = s^2.
\tag{4.11}
$$

Therefore, by Lyapunov's theorem (e.g., [2], p. 44), in order to show that

$$
\mathcal{L}(Z_{n,1}) \to N(0, s^2),
\tag{4.12}
$$

it suffices to prove that

$$
\frac{n}{\left(\sqrt{nh_n^{d+2}}\right)^4} E\left[ K\left(\frac{p - X}{h_n}\right) \langle f'(p), \mathcal{E}_p^{-1}(X) \rangle I(X \in \mathcal{B}_n) \right. \\
\left. - EK\left(\frac{p - X}{h_n}\right) \langle f'(p), \mathcal{E}_p^{-1}(X) \rangle I(X \in \mathcal{B}_n) \right]^4 \to 0
\tag{4.13}
$$

By the hypothesis on $h_n$ and (3.11) we can ignore the expected value within the square bracket, and for the rest, proceeding as usual, we have

$$\frac{1}{nh_n^{2(d+2)}} E\left[K\left(\frac{p-X}{h_n}\right)\langle f'(p), \mathcal{E}_p^{-1}(X)\rangle I(X \in \mathcal{B}_n)\right]^4$$
$$= \frac{1}{(4\pi)^2 |\mu| nh_n^{2(d+2)}} \int_{\tilde{\mathcal{B}}_n} e^{-(|v|^2+O(|v|^4))/h_n^2} \langle f'(p), v\rangle^4 (1+O(|v|^2)) dv$$
$$\leq \frac{2h_n^{d+4}}{(4\pi)^2 |\mu| nh_n^{2(d+2)}} \int_{\mathbf{R}^d} e^{-|v|^2/2} |f'(p)|^4 dv = O(1/(nh_n^d)) \to 0,$$

proving (4.13), and therefore, the limit (4.12). Now the theorem follows from Proposition 2.1, (4.5), (4.8), (4.9) and (4.12). □

This result extends without effort to the CLT uniform in $f$, which is the main result in this section.

**Theorem 4.2.** *Let $U$ be a normal neighborhood of $p$ and let $\mathcal{F}$ be a class of functions uniformly bounded up to the third order in U. Assume $nh_n^{d+4} \to 0$, and $nh_n^d \to \infty$. Then, as $n \to \infty$, the processes*

$$\text{(4.14)} \qquad \left\{\sqrt{nh_n^{d+2}} \Big[\Delta_{h_n,n} f(p) - \frac{1}{|\mu|}\Delta_M f(p)\Big] : \ f \in \mathcal{F}\right\}$$

*converge in law in $\ell^\infty(\mathcal{F})$ to the Gaussian process*

$$\text{(4.15)} \qquad \left\{G(f) := \frac{1}{2^{d/2}(2\pi)^{d/4}|\mu|^{1/2}} \sum_{j=1}^d Z_j \frac{\partial f}{\partial x_j^p}(p) : \ f \in \mathcal{F}\right\},$$

*where $Z = (Z_1, \ldots, Z_d)$ is the standard normal vector in $T_p(M)$ $(= \mathbf{R}^d)$.*

*Proof.* The proof of Proposition 4.1 applied to $f = \sum_{j=1}^r \alpha_j f_j$, with $f_j \in \mathcal{F}$, shows that the finite dimensional distributions of the processes (4.14) converge to those of the process (4.15) (by the definition (4.3) of $s = s(f)$). Also, by Theorem 3.1, we can center the processes (4.14). Hence, the Theorem will follow if we show that the processes $Z_n = Z_n(f)$ in (4.4) are asymptotically equicontinuous with respect to a totally bounded pseudometric on $\mathcal{F}$ (e.g., [5]).

First, by the computation (3.8) we can restrict the range of $X_i$ to $X_i \in \mathcal{B}_n$ by taking $L^2/4 > d+3$, because

$$\frac{1}{\sqrt{nh_n^{d+2}}} E \sup_{f\in\mathcal{F}} \bigg| \sum_{i=1}^n \bigg(K\bigg(\frac{p-X_i}{h_n}\bigg)(f(X_i)-f(p))I(X_i \in M\setminus\mathcal{B}_n)$$
$$- EK\bigg(\frac{p-X}{h_n}\bigg)(f(X_i)-f(p))I(X\in M\setminus\mathcal{B}_n)\bigg)\bigg|$$
$$\leq \frac{4cnh_n^{L^2/4}}{(4\pi)^{d/2}\sqrt{nh_n^{d+2}}} \to 0,$$

since $nh_n^{d+4} \to 0$.

Now that we can restrict to $X_i \in \mathcal{B}_n$, we only need to consider $Z_{n,i}(f)$, $i = 1, 2, 3$, as defined by equations (4.10), (4.6) and (4.7) respectively. Asymptotic equicontinuity of $Z_{n,1}(f)$ follows because

$$E \sup_{\|f'(p)\| \leq \delta} |Z_{n,1}(f)|^2 \leq \delta^2 E \left\| \frac{1}{\sqrt{n h_n^{d+2}}} \sum_{i=1}^n \left( K\left(\frac{p - X_i}{h_n}\right) \mathcal{E}_p^{-1}(X_i) I(X_i \in \mathcal{B}_n) \right. \right.$$
$$\left. \left. - EK\left(\frac{p-X}{h_n}\right) \mathcal{E}_p^{-1}(X) I(X \in \mathcal{B}_n) \right) \right\|^2$$
$$\leq \frac{\delta^2}{h_n^{d+2}} E \left\| K\left(\frac{p-X}{h_n}\right) \mathcal{E}_p^{-1}(X) I(X \in \mathcal{B}_n) \right\|^2$$
$$= \frac{\delta^2}{(4\pi)^d |\mu| h_n^{d+2}} \int_{\tilde{\mathcal{B}}_n} e^{-\|p - \mathcal{E}_p(v)\|^2 / 2h_n^2} |v|^2 \sqrt{\det(g_{ij})}(v) dv$$
$$\leq \frac{\delta^2}{(4\pi)^d |\mu|} \int_{\mathbf{R}^d} e^{-|v|^2/3} |v|^2 dv$$

which tends to zero when we take sup over $n$ and then limit as $\delta \to 0$.

Next, by the computation in (4.9),

$$E \sup_{f \in \mathcal{F}} |Z_{n,3}(f)| = O\left(\sqrt{n h_n^4}\right) \to 0.$$

Finally we consider $Z_{n,2}(f)$. Let $\varepsilon_i$ be i.i.d. Rademacher variables independent of $\{X_i\}$. Then, by symmetrization,

$$E \sup_{f \in \mathcal{F}} Z_{n,2}^2(f)$$
$$\leq \frac{4}{n h_n^{d+2}} E E_\varepsilon \sup_{f \in \mathcal{F}} \left| \sum_{i=1}^n \varepsilon_i K\left(\frac{p - X_i}{h_n}\right) f''(p)(\mathcal{E}_p^{-1}(X_i), \mathcal{E}_p^{-1}(X_i)) I(X_i \in \mathcal{B}_n) \right|^2.$$

Next, we recall that, for an operator $A$ in $\mathbf{R}^d$, (or in $T_p(M)$), we have the following identity for its quadratic form

$$A(u, v) := \langle Au, v \rangle = \langle A, u \otimes v \rangle_{HS},$$

where in orthonormal coordinates, $u \otimes v$ is the $d \times d$-matrix with entries $u_i v_j$, and the Hilbert-Schmidt inner product of two matrices is just the inner product in $\mathbf{R}^{2d}$. Also note in particular that $\|u \otimes v\|_{HS} = |u||v|$. Therefore,

$$E_\varepsilon \sup_{f \in \mathcal{F}} \left| \sum_{i=1}^n \varepsilon_i K\left(\frac{p - X_i}{h_n}\right) f''(p)(\mathcal{E}_p^{-1}(X_i), \mathcal{E}_p^{-1}(X_i)) I(X_i \in \mathcal{B}_n) \right|^2$$
$$= E_\varepsilon \sup_{f \in \mathcal{F}} \left\langle f''(p), \sum_{i=1}^n \varepsilon_i K\left(\frac{p - X_i}{h_n}\right) (\mathcal{E}_p^{-1}(X_i) \otimes \mathcal{E}_p^{-1}(X_i)) I(X_i \in \mathcal{B}_n) \right\rangle_{HS}^2$$
$$\leq E_\varepsilon \sup_{f \in \mathcal{F}} |f''(p)|_{HS}^2 E_\varepsilon$$
$$\quad \times \left\| \sum_{i=1}^n \varepsilon_i K\left(\frac{p - X_i}{h_n}\right) (\mathcal{E}_p^{-1}(X_i) \otimes \mathcal{E}_p^{-1}(X_i)) I(X_i \in \mathcal{B}_n) \right\|_{HS}^2$$
$$\leq b^2 \sum_{i=1}^n K^2\left(\frac{p - X_i}{h_n}\right) \|(\mathcal{E}_p^{-1}(X_i) \otimes \mathcal{E}_p^{-1}(X_i))\|_{HS}^2 I(X_i \in \mathcal{B}_n) =: c^2 \Lambda_n^2.$$

Now, by (4.8),

$$E\Lambda_n^2 = nEK^2\left(\frac{p-X}{h_n}\right)\|\mathcal{E}_p^{-1}(X)\|^4 I(X \in \mathcal{B}_n) \leq \frac{2nh_n^{d+4}}{(4\pi)^d|\mu|}\int_{\mathbf{R}^d} e^{-|v|^2/3}|v|^4 dv,$$

which gives

$$E \sup_{f \in \mathcal{F}} Z_{n,2}^2(f) = O(h_n^2) \to 0. \qquad \Box$$

A simpler proof along similar lines gives the following law of large numbers:

**Theorem 4.3.** *Let $U$ be a normal neighborhood of $p$ and let $\mathcal{F}$ be uniformly bounded and equicontinuous up to the second order in $U$. Assume $h_n \to 0$ and $nh_n^{d+2} \to \infty$. Then*

(4.16) $$\sup_{f \in \mathcal{F}} \left|\Delta_{h_n,n}f(p) - \frac{1}{|\mu|}\Delta_M f(p)\right| \to 0 \quad \text{in pr.}$$

A Law of the Iterated Logarithm is also possible, but we refrain from presenting one since in the next section we will give a law of the logarithm for the sup over $f \in \mathcal{F}$ and $p \in M$, and the same methods, with a simpler proof, give the LIL at a single point.

## 5. Uniform approximation of the Laplacian by graph Laplacians

This section is devoted to results about approximation of the Laplacian by graph Laplacians not only uniformly on the functions $f$, but also on the points $p \in M$, $M$ a compact submanifold or $M$ as in Remark 3.2. The distributional convergence requires extra work (recall the Bickel-Rosenblatt theorem on the asymptotic distribution of the sup of the difference between a density and its kernel estimator) and will not be considered here.

Although the results in this section are also valid in the situation of Remark 3.2, we will only state them for $M$ a compact submanifold (without boundary). Also, we will identify $M$ with $\phi(M)$, that is, the imbedding $\phi$ will not be displayed.

**Theorem 5.1.** *Let $M$ be a compact Riemannian submanifold of dimenison $d < m$ of $\mathbf{R}^m$, let $X, X_i$ be i.i.d. with law $\mu/|\mu|$ and let $K$ be as defined in (3.1). Let $\mathcal{F}$ be a class of functions uniformly bounded and equicontinuous up to the second order in $M$. If $h_n \to 0$ and $nh_n^{d+2}/\log h_n^{-1} \to \infty$, then*

(5.1) $$\sup_{f \in \mathcal{F}} \sup_{q \in M} \left|\Delta_{h_n,n}f(q) - \frac{1}{|\mu|}\Delta_M f(q)\right| \to 0 \quad \text{a.s.}$$

*as $n \to \infty$. Moreover, if $\mathcal{F}$ is a class of functions uniformly bounded up to the third order in $M$, and, in addition to the previous conditions on $h_n$, $nh_n^{d+4}/\log h_n^{-1} \to 0$, then*

(5.1′) $$\sup_{f \in \mathcal{F}} \sup_{q \in M} \left|\Delta_{h_n,n}f(q) - \frac{1}{|\mu|}\Delta_M f(q)\right| = O\left(\sqrt{\frac{\log(1/h_n)}{nh_n^{d+2}}}\right) \quad \text{as} \quad n \to \infty \quad \text{a.s.}$$

*Proof.* By Remark 3.1 on Theorem 3.1 (more precisely, by its uniform version), in order to prove (5.1) it suffices to show that

(5.2) $$\sup_{f \in \mathcal{F}} \sup_{q \in M} \left|\Delta_{h_n,n}f(q) - \Delta_{h_n}f(q)\right| \to 0 \quad a.s.$$

Let $\mathcal{B}_{n,q} = \{x \in M : \|x - q\| < Lh_n(\log h_n^{-1})^{1/2}\}$, where, we recall, $\|\cdot\|$ is the norm in $\mathbf{R}^m$. Then, as in (3.8) and (4.5), if $L^2/4 > d+2$,

$$
\begin{aligned}
\sup_f \sup_q &\left| \frac{1}{nh_n^{d+2}} \sum_{i=1}^n \Big( K((q-X_i)/h_n) I(X_i \in \mathcal{B}_{n,q}^c)(f(X_i) - f(q)) \right. \\
&\left. - E\Big( K((q-X)/h_n) I(X \in \mathcal{B}_{n,q}^c)(f(X) - f(q)) \Big) \Big) \right| \\
&\leq \frac{4\|f\|_\infty h_n^{L^2/4}}{h_n^{d+2}} \to 0.
\end{aligned}
\tag{5.3}
$$

To establish (5.1), we show that

$$
\begin{aligned}
E_n := \frac{1}{nh_n^{d+2}} E \sup_f \sup_q &\left| \sum_{i=1}^n \Big( K((q-X_i)/h_n) I(X_i \in \mathcal{B}_{n,q})(f(X_i) - f(q)) \right. \\
&\left. - EK((q-X)/h_n) I(X \in \mathcal{B}_{n,q})(f(X) - f(q)) \Big) \right| \to 0,
\end{aligned}
\tag{5.4}
$$

and use Talagrand's [14] concentration inequality to transform this into a statement on a.s. convergence.

Each function $f \in \mathcal{F}$ can be extended to a twice continuously differentiable function $\overline{f}$ on a compact domain $N$ of $\mathbf{R}^m$ with $M$ in its interior such that the classes $\{\overline{f} : f \in \mathcal{F}\}$ $\{\overline{f}' : f \in \mathcal{F}\}$ $\{\overline{f}'' : f \in \mathcal{F}\}$ are uniformly bounded and $\{\overline{f}'' : f \in \mathcal{F}\}$ is equicontinuous on $N$ (use a finite partition of unity to patch together convenient extensions of $f$ in each of the sets in a finite cover of $M$ by e.g., geodesic balls: see e.g. Lee [9], pp. 15-16). Then,

$$f(X_i) - f(q) = \overline{f}'(q + \theta(X_i - q))(X_i - q)$$

for some point $0 \leq \theta = \theta_{q,X_i} \leq 1$.

Note that, $M$ being compact, $\mathcal{B}_{n,q}$ is contained in one of a finite number of uniformly normal neighborhoods for all $n \geq N_0$, with $N_0 < \infty$ independent of $q$, so, we can use $q$-normal coordinates and notice that for these coordinates, on $\mathcal{B}_{n,q}$, we have the inequalities (3.10) holding uniformly in $q$ (by Proposition 2.2). Since the derivative $\overline{f}'$ is uniformly bounded, for $n \geq N_1$ (independent of $q$), we have

$$
\begin{aligned}
&EK^2((q-X_i)/h_n) I(X_i \in \mathcal{B}_{n,q})(f(X_i) - f(q))^2 \\
&\leq CEK^2((q-X_i)/h_n) I(X_i \in \mathcal{B}_{n,q}) \|X_i - q\|^2,
\end{aligned}
$$

which, in view of (3.10), can be further bounded by

$$
\begin{aligned}
&2CEK^2((q-X_i)/h_n) I(X_i \in \mathcal{B}_{n,q}) |\mathcal{E}_q^{-1}(X_i)|^2 \\
&\leq 2C \int_{\tilde{\mathcal{B}}_{n,q}} e^{-(|v|^2 - C|v|^4)/2h_n^2} |v|^2 (1 + C_1|v|^2) dv \\
&\leq 2C h_n^{d+2} \int_{\mathbf{R}^d} e^{-|v|^2/3} |v|^2 (1 + C_1|v|^2) dv \leq C_2 h_n^{d+2},
\end{aligned}
$$

so we have

$$EK^2((q-X_i)/h_n) I(X_i \in \mathcal{B}_{n,q})(f(X_i) - f(q))^2 \leq C_2 h_n^{d+2} \tag{5.5}$$

with a constant $C_2$ that does not depend on $q$.

To prove (5.4), we replace $f(X_i) - f(q)$ by its Taylor expansion of the second order:

$$f(X_i) - f(q) = \overline{f}'(q)(X_i - q) + \frac{1}{2}\overline{f}''(q)(X_i - q, X_i - q) + r_n(f; q; X_i),$$

where

$$\sup_{q \in M} \sup_{f \in \mathcal{F}} r_n(f; q; X) \leq \delta_n \|X - q\|^2$$

with $\delta_n \to 0$ as $n \to \infty$ (because of equicontinuity of $\{\overline{f}'' : f \in \mathcal{F}\}$ and the fact that $\|X - q\| < Lh_n (\log h_n^{-1})^{1/2} \to 0$).

The first order term leads to bounding the expectation

$$\frac{1}{nh_n^{d+2}} E \sup_f \sup_q \Big| \Big\langle \overline{f}'(q), \sum_{i=1}^n \Big( K((q - X_i)/h_n) I(X_i \in \mathcal{B}_{n,q})(X_i - q) \\ - EK((q - X)/h_n) I(X \in \mathcal{B}_{n,q})(X - q) \Big) \Big\rangle \Big|,$$

which is smaller than

$$\frac{b}{nh_n^{d+2}} E \sup_f \sup_q \Big\| \sum_{i=1}^n \Big( K((q - X_i)/h_n) I(X_i \in \mathcal{B}_{n,q})(X_i - q) \\ - EK((q - X)/h_n) I(X \in \mathcal{B}_{n,q})(X - q) \Big) \Big\|,$$

where $b$ is a uniform upper bound on $\overline{f}'$. Denote the coordinates of $x \in \mathbf{R}^m$ in the canonical basis of $\mathbf{R}^m$ by $x_\alpha$, $\alpha = 1, \ldots, m$ and consider the class of functions $M \mapsto \mathbf{R}$,

$$\mathcal{G} = \{f_{q,h,\lambda}(x) := e^{-\|q-x\|^2/4h^2} I(\|x - q\| < \lambda)(x_\alpha - q_\alpha) : q \in M, h > 0, \lambda > 0\}.$$

By arguments of Nolan and Pollard (1987), the class of functions of $x$, $\{e^{-\|q-x\|^2/4h^2} : q \in M, h > 0\}$ is VC subgraph; and it is well known that the open balls in $\mathbf{R}^m$ are VC and that the class of functions $\{x_\alpha - q_\alpha : q \in M\}$ is also VC subgraph (see, e.g., [5]). The three classes are bounded (resp. by 1, 1 and $2\sup\{\|x\| : x \in M\}$) and therefore, by simple bounds on covering numbers, the product of the three classes is VC-type with respect to the constant envelope $C = 2 + 2\sup\{\|x\| : x \in M\}$. In particular, if $N(\mathcal{G}, \varepsilon)$ are the covering numbers for $\mathcal{G}$ in $L_2$ of any probability measure, then

$$N(\mathcal{G}, \varepsilon) \leq \left(\frac{A}{\varepsilon}\right)^v$$

for some $A, v < \infty$ and all $\varepsilon$ less than or equal to the diameter of $\mathcal{G}$. Hence, by inequality (2.2) in [7], there exists a constant $R$ such that

$$\begin{aligned}
& E \sup_q \Big| \sum_{i=1}^n \Big( K((q - X_i)/h_n) I(\|X_i - q\| < Lh_n(\log h_n^{-1})^{1/2})(X_{i,\alpha} - q_\alpha) \\
& \quad - E\Big( K((q - X)/h_n) I(\|X - q\| < Lh_n(\log h_n^{-1})^{1/2})(X_\alpha - q_\alpha) \Big| \\
& \leq R\Big(\sqrt{n}\sigma\sqrt{\log \frac{A}{\sigma}} \vee \log \frac{A}{\sigma}\Big),
\end{aligned}$$

(5.6)

where $\sigma^2 \geq \sup_{f \in \mathcal{G}} Ef^2(X)$. Now, to compute $\sigma$ we use again our observations before the proof of (5.5). For $n \geq N_1$ (independent of $q$), we have

$$\sup_{f \in \mathcal{G}} Ef^2(X) \leq \int_{\tilde{\mathcal{B}}_{n,q}} e^{-(|v|^2 - C|v|^4)/2h_n^2} |v|^2 (1 + C_1|v|^2) dv$$

$$\leq h_n^{d+2} \int_{\mathbf{R}^d} e^{-|v|^2/3} |v|^2 dv \leq C_2 h_n^{d+2},$$

for some $C_2 < \infty$ independent of $q$. So, we can take $\sigma^2 = C_2 h_n^{d+2}$. Hence, by the hypothesis on $h_n$, the right hand side of (5.6) is bounded by

$$R'\Big(nh_n^{d+2} \log \frac{A}{h_n}\Big)^{1/2}$$

for some $R' < \infty$.

To handle the second order term, note that

$$\frac{1}{nh_n^{d+2}} E \sup_f \sup_q \Big| \sum_{i=1}^n \Big(K((q-X_i)/h_n)I(X_i \in \mathcal{B}_{n,q})\overline{f}''(q)(X_i - q, X_i - q)$$

$$- EK((q-X)/h_n)I(X \in \mathcal{B}_{n,q})\overline{f}''(q)(X - q, X - q)\Big)\Big|$$

$$= \frac{1}{nh_n^{d+2}} E \sup_f \sup_q \Big| \Big\langle \overline{f}''(q),$$

$$\sum_{i=1}^n \Big(K((q-X_i)/h_n)I(X_i \in \mathcal{B}_{n,q})(X_i - q) \otimes (X_i - q)$$

$$- EK((q-X)/h_n)I(X \in \mathcal{B}_{n,q})(X - q) \otimes (X - q)\Big)\Big\rangle_{HS}\Big|,$$

which is dominated by

$$\frac{b}{nh_n^{d+2}} E \sup_f \sup_q \Big\| \sum_{i=1}^n \Big(K((q-X_i)/h_n)I(X_i \in \mathcal{B}_{n,q})(X_i - q) \otimes (X_i - q)$$

$$- EK((q-X)/h_n)I(X \in \mathcal{B}_{n,q})(X - q) \otimes (X - q)\Big)\Big\|_{HS}$$

(with $b$ being a uniform upper bound on $\overline{f}''$). Here $\otimes$ denotes the tensor product of vectors of $\mathbf{R}^m$ and $\|\cdot\|_{HS}$ is the Hilbert-Schmidt norm for linear transformations of $\mathbf{R}^m$. This leads to bounding the expectation

$$E \sup_q \Big| \sum_{i=1}^n \Big(K((q-X_i)/h_n)I(\|X_i - q\| < Lh_n(\log h_n^{-1})^{1/2})(X_{i,\alpha} - q_\alpha)(X_{i,\beta} - q_\beta)$$

(5.6′)
$$- E\Big(K((q-X)/h_n)I(\|X - q\| < Lh_n(\log h_n^{-1})^{1/2})(X_\alpha - q_\alpha)(X_\beta - q_\beta)\Big)\Big|$$

for all $1 \leq \alpha, \beta \leq m$, which is done using the inequality for empirical processes on VC-subgraph classes exactly the same way as in the case (5.6). This time the bound becomes

$$R'\Big(nh_n^{d+4} \log \frac{A}{h_n}\Big)^{1/2}.$$

For the remainder, we have the bound

$$\frac{1}{nh_n^{d+2}} E \sup_f \sup_q \Big| \sum_{i=1}^n K((q-X_i)/h_n)I(X_i \in \mathcal{B}_{n,q})r_n(f,q,X_i)$$
$$- EK((q-X)/h_n)I(X \in \mathcal{B}_{n,q})r_n(f,q,X)\Big|$$
$$\leq \frac{\delta_n}{nh_n^{d+2}} E \sup_q \sum_{i=1}^n \Big( K((q-X_i)/h_n)I(X_i \in \mathcal{B}_{n,q})\|X_i-q\|^2$$
$$+ n\frac{\delta_n}{nh_n^{d+2}} \sup_q E\Big( K((q-X)/h_n)I(X \in \mathcal{B}_{n,q})\|X-q\|^2 \Big)$$
$$\leq \frac{\delta_n}{nh_n^{d+2}} E \sup_q \sum_{i=1}^n \Big( K((q-X_i)/h_n)I(X_i \in \mathcal{B}_{n,q})\|X_i-q\|^2$$
$$- EK((q-X)/h_n)I(X \in \mathcal{B}_{n,q})\|X-q\|^2 \Big)$$
$$+ \frac{2\delta_n}{h_n^{d+2}} \sup_q EK((q-X)/h_n)I(X \in \mathcal{B}_{n,q})\|X-q\|^2.$$

The first expectation is bounded again by using the inequality for VC-subgraph classes and the bound in this case is

$$\frac{\delta_n}{nh_n^{d+2}} \Big( nh_n^{d+4} \log \frac{A}{h_n} \Big)^{1/2} = \delta_n \sqrt{\frac{\log(A/h_n)}{nh_n^d}} \to 0.$$

The second expectation is bounded by replacing $\|X-q\|^2$ by $|\mathcal{E}_q^{-1}(X)|^2$ and changing variables in the integral (as it has been done before several times). This yields a bound of the order $C\delta_n$, which also tends to 0.

Combining the above bounds establishes (5.4). One of the versions of Talagrand's inequality (e.g., [10]) together with (5.5) gives with some constant $K > 0$ and with probability at least $1 - e^{-t}$

$$\frac{1}{n} \sup_f \sup_q \Big| \sum_{i=1}^n \Big( K((q-X_i)/h_n)I(X_i \in \mathcal{B}_{n,q})(f(X_i)-f(q))$$
$$- EK((q-X)/h_n)I(X \in \mathcal{B}_{n,q})(f(X)-f(q))\Big)\Big|$$
$$\leq K\Big( h_n^{d+2} E_n + \sqrt{h_n^{d+2}\frac{t}{n}} + \frac{t}{n} \Big).$$

Taking $t := t_n := A \log n$ with large enough $A$, so that $\sum_n e^{-t_n} < \infty$ and using Borel-Cantelli Lemma shows that a.s. for large enough $n$

$$\frac{1}{nh_n^{d+2}} \sup_f \sup_q \Big| \sum_{i=1}^n \Big( K((q-X_i)/h_n)I(X_i \in \mathcal{B}_{n,q})(f(X_i)-f(q))$$
$$- EK((q-X)/h_n)I(X \in \mathcal{B}_{n,q})(f(X)-f(q))\Big)\Big|$$
$$\leq K\Big( E_n + \sqrt{\frac{A\log n}{nh_n^{d+2}}} + \frac{A\log n}{nh_n^{d+2}} \Big)$$

and since, in view of (5.4) and under the condition $nh_n^{d+2}/\log h_n^{-1} \to \infty$, the right hand side tends to 0. This and (5.3) yield (5.1). (Note that $nh_n^{d+2}/\log h_n^{-1} \to \infty$ implies $nh_n^{d+2}/\log n \to \infty$.)

The proof of (5.1') requires the following version of Taylor's expansion of $\overline{f}$:

$$f(X_i) - f(q) = \overline{f}'(q)(X_i - q) + \frac{1}{2}\overline{f}''(q)(X_i - q, X_i - q)$$
$$+ \frac{1}{6}\overline{f}'''(q + \theta_i(X_i - q))(X_i - q, X_i - q, X_i - q).$$

The first two terms have been handled before, and the expectations of the sup-norms of the corresponding empirical processes were shown to be $O(\sqrt{\frac{\log h_n^{-1}}{nh_n^{d+2}}})$. The third term leads to bounding

$$\frac{1}{nh_n^{d+2}} E \sup_f \sup_q \Big| \sum_{i=1}^n \Big( K((q - X_i)/h_n) I(X_i \in \mathcal{B}_{n,q})$$
$$- EK((q - X)/h_n) I(X \in \mathcal{B}_{n,q}) \Big)$$
$$\times \frac{1}{6}\overline{f}'''(q + \theta_i(X_i - q))(X_i - q, X_i - q, X_i - q) \Big|,$$

which, for $f'''$ uniformly bounded by $b$, is smaller than

$$\frac{b}{6nh_n^{d+2}} E \sup_q \sum_{i=1}^n \Big( K((q - X_i)/h_n) I(X_i \in \mathcal{B}_{n,q}) \|X_i - q\|^3$$
$$+ \frac{bn}{6nh_n^{d+2}} \sup_q EK((q - X)/h_n) I(X \in \mathcal{B}_{n,q}) \|X - q\|^3 \Big)$$
$$\leq \frac{b}{6nh_n^{d+2}} E \sup_q \Big| \sum_{i=1}^n \Big( K((q - X_i)/h_n) I(X_i \in \mathcal{B}_{n,q}) \|X_i - q\|^3$$
$$- EK((q - X)/h_n) I(X \in \mathcal{B}_{n,q}) \|X - q\|^3 \Big) \Big|$$
$$+ \frac{b}{3h_n^{d+2}} \sup_q EK((q - X)/h_n) I(X \in \mathcal{B}_{n,q}) \|X - q\|^3,$$

which can be handled exactly as before and shown to be of the order

$$h_n \sqrt{\frac{\log h_n^{-1}}{nh_n^2}} + h_n = o\Big(\sqrt{\frac{\log h_n^{-1}}{nh_n^{d+2}}}\Big),$$

by the conditions on $h_n$. Using Talagrand's inequality the same way as before, completes the proof of (5.1'). □

We conclude with the following theorem, whose proof is a little longer and more involved, but it is based on a methodology that is well known and well described in the literature (see [6] and [7]). Its extension to the case of manifolds requires some work, but is rather straightforward.

**Theorem 5.2.** *Let $M$ be a compact Riemannian submanifold of dimenison $d < m$ of $\mathbf{R}^m$, let $X, X_i$ be i.i.d. with law $\mu/|\mu|$ and let $K$ be as defined in (3.1). Assume that $h_n \to 0$, $nh_n^{d+2}/\log h_n^{-1} \to \infty$, and $nh_n^{d+4}/\log h_n^{-1} \to 0$. Let $\mathcal{F}$ be a class of* *functions uniformly bounded up to the third order in M. Then,*

$$\lim_{n\to\infty} \sqrt{\frac{nh_n^{d+2}}{2\log h_n^{-d}}} \sup_{f\in\mathcal{F}} \sup_{q\in M} \left|\Delta_{h_n,n}f(q) - \frac{1}{|\mu|}\Delta_M f(q)\right|$$
$$= \frac{\sup_{f\in\mathcal{F}, q\in M}\left(\sum_{j=1}^d \left(\frac{\partial f}{\partial x_j^q}(q)\right)^2\right)^{1/2}}{2^{d/2}(2\pi)^{d/4}|\mu|^{1/2}} \quad \text{a.s.,}$$

*where $x_j^q$ denote normal coordinates centered at q.*